\documentclass[12pt]{amsart}
\usepackage{amssymb,times}
\newtheorem{thm}{Theorem}
\newtheorem{prop}[thm]{Proposition}
\newtheorem{lem}[thm]{Lemma}
\newtheorem{cor}[thm]{Corollary}

\theoremstyle{remark}
\newtheorem{rem}[thm]{Remark}
\newtheorem{ex}[thm]{Example}

\theoremstyle{definition}
\newtheorem{defn}[thm]{Definition}

\newcommand{\C}{\mathbb{ C}}

\newcommand{\R}{\mathbb{ R}}

\newcommand{\Z}{\mathbb{ Z}}

\newcommand{\Spc}{\mathrm{Spin}^c}

\newcommand{\Ha}{\mathcal{ H}}

\newcommand{\fs}{\mathfrak s}

\newcommand{\DA}{{D_{\! A}^+}}

\title{Monopole classes and Einstein metrics}
\author{D.~Kotschick}
\address{Mathematisches Institut, Ludwig-Maximilians-Universit\"at M\"unchen,
Theresienstr.~39, 80333 M\"unchen, Germany}
\email{dieter@member.ams.org}
\email{Jan.Wehrheim@mathematik.uni-muenchen.de}
\thanks{I am grateful to V.~Braungardt, M.~Furuta, P.~Kronheimer and 
J.~Wehrheim for useful conversations. The author is a member of the 
{\sl European Differential Geometry Endeavour} (EDGE), Research Training 
Network HPRN-CT-2000-00101, supported by The European Human Potential Programme}
\date{March 19, 2003, revised May 3, 2003; MSC 2000 classification: primary 57R57; secondary 53C25, 57R55}

\begin{document}

\begin{abstract}
We introduce the notion of a special monopole class on a four-manifold. 
This is used to prove restrictions on the smooth structures of Einstein 
manifolds. As an application we prove that there are Einstein four-manifolds 
which are simply connected, spin, and satisfy the strict Hitchin--Thorpe 
inequality, and which are homeomorphic to manifolds without Einstein metrics. 
\end{abstract}

\maketitle

\section{Introduction}

Seiberg--Witten theory has been developed as a theory of moduli spaces 
of solutions to the monopole equations by imitating the development for 
the anti-self-dual Yang--Mills equations. Therefore, in the applications 
the use of the whole moduli space has been emphasized, following the 
ideas of Donaldson in Yang--Mills theory. Nevertheless, unlike Donaldson 
theory, Seiberg--Witten theory has incorporated geometric arguments 
exploiting a single solution to the monopole equations for explicit 
constructions, without necessarily studing all its deformations. Rather 
than imitating Donaldson theory, these arguments are more reminiscent of 
Yau's work on the Calabi conjectures, where a single Einstein metric is 
constructed and used as input in geometric arguments, or of the use of 
minimal surfaces in three-dimensional topology or in the theory of manifolds 
of positive scalar curvature. Going back even further, there is a clear 
parallel with the use of harmonic forms and harmonic spinors in geometry. 
Of course, those are solutions of linear partial differential equations, 
rather than non-linear ones, and one might speculate that the mildness of 
the non-linearity of the Seiberg--Witten equations manifests itself in 
the way the solutions are used\footnote{With the hindsight of 
Seiberg--Witten theory, it is worth rereading early appraisals of 
Donaldson theory, such as~\cite{At}.}.

For arguments using the existence of a solution to the monopole 
equations rather than the non-triviality of a topological invariant, 
Kronheimer~\cite{Kr} introduced the notion of a monopole class. This 
is the characteristic class of a $\Spc$-structure for which there are 
solutions of the monopole equations for arbitrary metrics.
The usefulness of this notion stems from the fact that 
one can make geometric arguments with a solution, without ever studying 
its deformation theory. The purpose of this paper is to derive certain 
constraints on the smooth structures underlying Einstein manifolds from 
the existence of monopole classes. In some of our arguments we need two 
additional properties of monopole classes, which are always satisfied by 
the basic classes coming from non-trivial Seiberg--Witten invariants, but 
which are not part of Kronheimer's definition of monopole classes. Thus, 
in Section~\ref{s:mono} we discuss variations on the definition, and the 
relations between these different definitions. This leads to issues of 
genericity and transversality for the monopole equations. 

In Section~\ref{s:Einstein} we derive restrictions on Einstein 
four-manifolds with monopole classes, and use these to give various 
examples of simply connected Einstein manifolds homeomorphic to smooth 
manifolds which cannot admit any Einstein metric. 

The existence of such examples is not new. I gave the first examples 
in~\cite{K}, showing that an obstruction to the existence of Einstein 
metrics discovered by LeBrun~\cite{lebrun} depends on the smooth structure, 
rather than being homotopy-invariant. Recently, Ishida--LeBrun~\cite{IL} 
asked if there are such examples which are spin. They found an obstruction 
that is applicable to spin manifolds, but were unable to show that it is 
not homotopy-invariant\footnote{See Remark~\ref{r:homotopy} below for an 
instance where it does turn out to be homotopy-invariant.}. We shall resolve 
their conundrum by exhibiting a simply connected Einstein manifold which is 
spin and satisfies the strict Hitchin--Thorpe inequality, and is homeomorphic 
to manifolds not admitting any Einstein metric because of the obstruction 
discussed in~\cite{IL}. 

Of course, the spin condition is a red herring. It is clear from the known 
examples that the obstructions to Einstein metrics derived from the 
Seiberg--Witten monopole equations have to do with the reducibility of the 
underlying manifolds, spin or non-spin. (Compare~\cite{lebrun,K,lebrun3,IL,BK} 
and Section~\ref{s:Einstein} below.) In fact, it is 
quite tempting to conjecture that, except for blowups of the complex 
projective plane, Einstein four-manifolds should be irreducible. In 
Section~\ref{s:irred} we prove some partial results in this direction.

In the Appendix, written jointly with J.~Wehrheim, we prove some results 
about Ricci-flat four-manifolds motivated by the conjecture that all closed 
Ricci-flat four-manifolds should be covered by $K3$ or $T^{4}$. From these 
results we conclude that the existence of strongly scalar-flat four-manifolds 
which are not covered by $T^{4}$ or $K3$ would imply the existence of some 
new obstruction to positive scalar curvature.


\section{Monopole classes}\label{s:mono}

Consider a closed smooth oriented four-manifold $X$ with a $\Spc$-structure 
$\fs$. For every choice of Riemannian metric $g$, the Seiberg--Witten 
monopole equations for $(X,\fs)$ with respect to $g$ are a system of coupled 
equations for a pair $(A,\Phi)$, where $A$ is a $\Spc$-connection in the spin 
bundle for $\fs$ and $\Phi$ is a section of the positive spin bundle $V_{+}$.
The equations are:
\begin{equation}\label{eq:D}
     \DA \Phi =0 \ ,
\end{equation}
\begin{equation}\label{eq:curv}
     F_{\hat A}^+ =\sigma (\Phi ,\Phi ) \ ,
\end{equation}
with $\DA$ the half-Dirac operator defined on spinors of positive chirality, 
and ${\hat A}$ the connection induced by $A$ on the determinant of the spin 
bundle. The right-hand side of the curvature equation~\eqref{eq:curv} is the 
$2$-form which, under the Clifford module structure determined by $\fs$, 
corresponds to the trace-free part of $\Phi\otimes\Phi^{*}$.

Solutions $(A,\Phi)$ with $\Phi\equiv 0$ are called {\it reducible}. If 
there is a reducible solution, then $c=c_{1}(\fs)$ is represented by an 
anti-self-dual harmonic form because of~\eqref{eq:curv}. This implies 
$c^{2}\leq 0$, with equality if and only if $c$ is a torsion class.

The following definition is due to Kronheimer~\cite{Kr}.
\begin{defn}
    A class $c\in H^{2}(X,\Z)$ is called a {\it monopole class}, if there 
    is a $\Spc$-structure $\fs$ on $X$ with $c=c_{1}(\fs)$ for which the 
    monopole equations~\eqref{eq:D} and~\eqref{eq:curv} admit a solution 
    $(A,\Phi)$ for every choice of metric $g$. 
    \end{defn}
Of course, on manifolds for which the Seiberg--Witten invariants are 
well-defined, every basic class is a monopole class. The rationale for 
considering the concept of a monopole class is that the existence of 
solutions to the monopole equations has immediate consequences, even when 
the corresponding invariants vanish. For example, Kronheimer~\cite{Kr} 
showed that the adjunction inequality $2g(\Sigma)-2\geq c\cdot\Sigma$ holds 
for any monopole class $c$ and any smoothly embedded surface $\Sigma\subset X$ 
of positive genus $g(\Sigma)$ which has trivial normal bundle.

Another immediate consequence, going back to Witten~\cite{Wi}, is that a 
non-torsion monopole class $c$ rules out the existence of a metric of positive 
scalar curvature as soon as $b_{2}^{+}(X)>0$, i.~e.~the intersection 
form of $X$ is not negative definite. This is seen by combining the 
Weitzenb\"ock formula for the Dirac operator with the curvature equation. 
A short calculation then shows that either $\Phi\equiv 0$, or one has the 
estimate $\vert\Phi(p)\vert\leq -s_{g}(p)$ at points $p\in X$ where $\vert\Phi\vert$ 
has a local maximum. Therefore, if the scalar curvature $s_{g}$ is non-negative, 
every solution of the monopole equations is reducible. If $c^{2}>0$, this 
is a contradiction. If $c^{2}\leq 0$, one has to perturb the metric, 
using that positivity of the scalar curvature is an open condition in the 
space of metrics. The following lemma then concludes the argument. A proof can 
be found in~\cite{DK}.
\begin{lem}[Donaldson]\label{l:D}
    Fix a non-torsion class $c\in H^{2}(X,\Z)$. If $b_{2}^{+}(X)>0$, then for a 
    generic metric $g$, there is no anti-self-dual harmonic form 
    representing the image of $c$ in $H^{2}(X,\R)$.
    \end{lem}

    Note that the above argument rules out metrics of zero scalar curvature 
if the monopole class $c$ satisfies $c^{2}>0$, but not otherwise, because 
scalar-flat metrics do not form an open set in the space of metrics. 


There are two fundamental properties of basic classes, the use of which has 
already become second nature to workers in the field, and which it would be 
desirable to have at hand for monopole classes. One is the non-negativity 
of the formal dimension of the associated moduli space of solutions to the 
monopole equations, the other is finiteness of the set of basic, 
respectively monopole classes. Therefore, we make the following definition, 
modifying Kronheimer's~\cite{Kr}, to ensure these extra properties.
\begin{defn}
    A {\it special monopole class} $c$ is a monopole class for which 
    $c^{2}\geq 2e(X)+3\sigma(X)$.
    \end{defn}
Here $e$ and $\sigma$ stand for the Euler characteristic and the 
signature. The inequality $c^{2}\geq 2e(X)+3\sigma(X)$ is equivalent 
to the non-negativity of the formal or expected dimension 
$$
d=\frac{1}{4}\left(c^{2}-(2e(X)+3\sigma(X))\right)
$$
of the solution space up to gauge equivalence. As an immediate 
consequence of the definition we have finiteness of special monopole 
classes:
\begin{lem}[Witten~\cite{Wi}]\label{l:finite}
    On any smooth closed oriented four-manifold $X$ there are at most 
    finitely many special monopole classes.
    \end{lem}
\begin{proof}
    Suppose $c$ is a special monopole class, and $(A,\Phi)$ is any 
    solution for the corresponding $\Spc$-structure $\fs$ and some 
    Riemannian metric $g$.
    The estimate on solutions we discussed before implies
\begin{equation}\label{eq:W}
    \vert\vert F_{\hat A}^+\vert\vert^{2}_{L^2} \leq \frac{1}{8}\vert\vert
s_g\vert\vert^{2}_{L^2} \ .
\end{equation}
The Chern--Weil formula 
$$
\vert\vert F_{\hat A}^-\vert\vert_{L^2}^2=\vert\vert F_{\hat A}^+
\vert\vert_{L^2}^2-4\pi^2c^2
$$ 
shows that an $L^2$ bound on $F_{\hat A}^+$, such as~\eqref{eq:W}, together 
with the assumption that $c^{2}$ is bounded below, implies an $L^2$ bound for 
$F_{\hat A}^-$.

Thus the projections of $F_{\hat A}$ to the harmonic subspaces
$\Ha^2_+$ and $\Ha^2_-$ are both bounded, so that the image
of $c$ in $H^2(X,\R )=\Ha^2_+\oplus\Ha^2_-$ is 
contained in a bounded set. Thus there are only finitely
many possibilities for $c$.
\end{proof}
Note that the proof works with any lower bound on $c^{2}$, 
equivalently any lower bound on the formal dimension $d$. 

\begin{rem}\label{r:IL}
    Ishida--LeBrun~\cite{IL} (p.~232) wrongly assert that Witten 
    proved finiteness for {\it all} monopole classes, not necessarily 
    special ones. Witten's argument, see~\cite{Wi} (p.~782/783), is 
    exactly the one in the proof of Lemma~\ref{l:finite} above, with 
    the additional assumption that $d=0$, cf.~(3.6) at the top of page 
    783, because he was working with the basic classes for 
    Seiberg--Witten invariants coming from zero-dimensional moduli 
    spaces. 
    \end{rem}
Here is an example with an infinite number of monopole classes, 
showing that not all monopole classes are special.
\begin{rem}
    On $X=\overline{\C P^{2}}$, every odd multiple of the generator 
    of $H^{2}(X,\Z)$ is a monopole class. For metrics of positive 
    scalar curvature on $X$, like the Fubini--Study metric, all 
    solutions are reducible.
    \end{rem}
In fact, this remark applies to any manifold with negative definite 
intersection form. Using the solution spaces of the monopole equations, 
one can prove Donaldson's theorem concluding that these intersection forms 
are diagonal over the integers, as soon as one has special monopole classes 
on these manifolds. That special monopole classes would have to exist for 
every negative definite four-manifold with non-standard intersection form 
is the content of Elkies's theorem~\cite{el}.

One should keep in mind that on any spin four-manifold $c=0$ is a 
monopole class, because for every metric one has the trivial reducible 
solution for a $\Spc$-structure induced by a spin structure. Considering 
connected sums of $S^{2}\times S^{2}$ and $S^{3}\times S^{1}$, one 
sees that the zero class does not obstruct the existence of positive 
scalar curvature metrics, and that the expected dimension of the 
corresponding moduli space can be positive, zero, or negative.
Monopole classes which are torsion play a different role from 
non-torsion classes.

As soon as $b_{2}^{+}(X)>0$, a non-torsion monopole class guarantees 
the existence of only irreducible solutions for generic metrics because 
of Lemma~\ref{l:D}. However, it is not known whether these are cut out 
transversely by the monopole equations. This is the issue of the 
missing ``generic metrics theorem'', which would be the analog of the 
Freed--Uhlenbeck theorem~\cite{FU} for the anti--self--dual Yang--Mills 
equations. If such a result were true, saying that for irreducible 
solutions of~\eqref{eq:D} and~\eqref{eq:curv} for a generic metric $g$ 
and a non-torsion $c=c_{1}(\fs)$ the linearized equation has no cokernel, 
then it would follow that on manifolds with $b_{2}^{+}(X)>0$ all 
non-torsion monopole classes are special monopole classes. 

Another useful variation of the definition is the following:
\begin{defn}
    A {\it generic monopole class} $c$ is a monopole class for which 
    there are solutions of the equations 
    \begin{equation}\label{eq:DD}
     \DA \Phi =0 \ ,
\end{equation}
\begin{equation}\label{eq:curv2}
     F_{\hat A}^+ =\sigma (\Phi ,\Phi )+\omega^{+} \ 
\end{equation}
for all metrics $g$ and all two-forms $\omega$.
    \end{defn}

\begin{lem}
    If $b_{2}^{+}(X)>0$, then generic monopole classes are special.
    \end{lem}
    \begin{proof}
	The usual transversality proof~\cite{thom,M} shows that for a 
	generic pair $(g,\omega)$, all the irreducible solutions of the 
	perturbed equations~\eqref{eq:DD} and~\eqref{eq:curv2} are cut out 
	transversely. Thus $d\geq 0$ as soon as there are irreducible 
	solutions for generic parameters. In view of $b_{2}^{+}(X)>0$, 
	irreducible solutions exist for any non-torsion monopole class by 
	Lemma~\ref{l:D}. But even if $c$ is a torsion class, the assumption 
	$b_{2}^{+}(X)>0$ allows us to choose generic parameters $(g,\omega)$ 
	for which there are no reducible solutions.
\end{proof}
If one does not have an invariant constructed from the moduli space 
which is unchanged under the perturbation by $\omega$, it seems quite 
hard to determine whether a given special monopole class is a generic 
one. Thus the monopole classes detected by three-dimensional topology that 
Kronheimer~\cite{Kr} discussed are all special, but whether they are generic 
is not clear to me. Recall that in many cases, the zero class is a 
monopole class, but is not a generic monopole class, because the 
trivial solution disappears when we perturb the curvature equation 
with a non-zero $\omega$. However, if a torsion class 
is a generic monopole class, so that not all solutions disappear 
under perturbation of the equations with a small non-zero $\omega$, 
then this monopole class does obstruct the existence of positive 
scalar curvature metrics if $b_{2}^{+}>0$, cf.~Lemma~\ref{l:D}.
Thus, rephrasing Witten's discussion in~\cite{Wi} in the language of 
monopole classes, we have:
\begin{prop}[Witten~\cite{Wi}]\label{p:Wi}
    Let $X$ be a closed oriented four-manifold with $b_{2}^{+}(X)>0$. 
    If $X$ admits a monopole class which is either non-torsion or 
    generic, then $X$ does not support any metric of positive scalar 
    curvature.
    \end{prop}
    
Obviously basic classes detected by numerical Seiberg--Witten invariants 
are generic monopole classes. A larger class of monopole classes, many 
of which are not basic, is detected by the stable cohomotopy invariant 
of Bauer--Furuta~\cite{BF}. These are generic, and therefore special, and 
often have $c^{2}>2e(X)+3\sigma(X)$, equivalently the formal dimensions 
$d$ of the moduli spaces are strictly positive. The following theorem 
summarizes and paraphrases some consequences of the connected sum formula 
for the stable cohomotopy invariant.
\begin{thm}[Bauer~\cite{B}]\label{t:B}
    Let $X_{i}$ be smooth closed oriented four-manifolds with 
    $b_{1}=0$ and $b_{2}^{+}\equiv 3\pmod 4$. Assume that $c_{i}$ is a 
    basic class on $X_{i}$ giving rise to a zero-dimensional moduli 
    space with odd numerical Seiberg--Witten invariant. Then 
    $c_{1}+c_{2}$, respectively $c_{1}+c_{2}+c_{3}$, is a generic 
    monopole class on $X_{1}\# X_{2}$, respectively on $X_{1}\# X_{2}\# 
    X_{3}$. 
    
    Similarly, $c_{1}+\ldots + c_{4}$ is a generic monopole class on 
    $X=X_{1}\#\ldots\# X_{4}$ if $b_{2}^{+}(X)\equiv 4\pmod 8$.
    \end{thm}
The formal dimensions of the relevant moduli spaces are $1$, $2$ and 
$3$ respectively.


\section{The exotic smooth structures of Einstein manifolds}\label{s:Einstein}

Recall that on simply connected four-manifolds the only classical 
obstruction to the existence of Einstein metrics is the Hitchin--Thorpe 
inequality
\begin{equation}\label{eq:HT}
e (X) \geq \frac{3}{2}\vert\sigma (X)\vert \ .
\end{equation}
Hitchin~\cite{HT} showed that Einstein manifolds for which~\eqref{eq:HT} 
is an equality are either flat or are isometric quotients of a $K3$ surface 
with a Calabi--Yau metric. Thus, the existence of smooth manifolds homeomorphic 
but not diffeomorphic to the $K3$ surface provides examples of homeomorphic 
manifolds such that one admits an Einstein metric and the other one does not. 
I gave the first such examples satisfying the strict Hitchin--Thorpe inequality 
in~\cite{K}. We shall paraphrase the argument to exclude Einstein metrics in 
certain cases in the next section. Suffice it to say here that it uses Witten's 
estimate~\eqref{eq:W} which implies that for any monopole class $c$ on $X$ and 
any metric $g$ the projection $c^{+}$ of $c$ to the self-dual summand in 
$H^2(X,\R )=\Ha^2_+\oplus\Ha^2_-$ satisfies 
$$
	(c^{+})^{2}\leq \frac{1}{32\pi^2}\vert\vert s_g\vert\vert^{2}_{L^2} \ .
$$
Sometimes the following can be used in the same way to obtain 
stronger results.
\begin{thm}[LeBrun~\cite{lebrun3}]\label{t:lebrun}
    Let $c$ be a monopole class on $X$. Then for any metric $g$ the 
    projection $c^{+}$ of $c$ to the self-dual summand in 
    $H^2(X,\R )=\Ha^2_+\oplus\Ha^2_-$ satisfies 
    \begin{equation}\label{eq:sharp}
	\frac{2}{3}(c^{+})^{2}\leq \frac{1}{4\pi^2}\left(\frac{1}{24}\vert\vert 
	s_g\vert\vert^{2}_{L^2}+2\vert\vert W_{+}\vert\vert^2_{L^2} \right) \ .
	\end{equation}
	In the case of equality the metric is almost K\"ahler.
    \end{thm}
As an immediate consequence, we have:
\begin{thm}\label{t:main}
    Let $c$ be a monopole class on $X$, and $d$ the expected 
    dimension of the corresponding moduli space. If $X$ admits an 
    Einstein metric, then $2e(X)+3\sigma(X)\geq 8d$. 
    In the case of equality $X$ admits a symplectic structure.
    \end{thm}
\begin{proof}
    If $g$ is an Einstein metric on $X$, the Gauss--Bonnet/Chern--Weil 
    formulae for the Euler characteristic and the signature combine to 
    give: 
    $$
    2e(X)+3\sigma(X)=\frac{1}{4\pi^2}\left(\frac{1}{24}\vert\vert 
	s_g\vert\vert^{2}_{L^2}+2\vert\vert W_{+}\vert\vert^2_{L^2} \right) \ .
    $$
    Combining this with~\eqref{eq:sharp} and
    $(c^{+})^{2}\geq c^{2}=4d+2e(X)+3\sigma(X)$ gives the result.
    
    The statement about the case of equality follows from the 
    corresponding statement in Theorem~\ref{t:lebrun}.
    \end{proof}
Note that whenever $d\leq 0$, this is weaker than the Hitchin--Thorpe 
inequality~\eqref{eq:HT}. In particular, the result is empty for 
monopole classes which are not special. 

Theorem~\ref{t:main} clarifies the arguments of Ishida--LeBrun~\cite{IL}. 
It could have appeared in their paper, but did not, because they mixed up 
the consequences of the existence of a monopole class with the device used 
to find such classes. One of their statements is the following:
\begin{cor}[Ishida--LeBrun~\cite{IL}]\label{c:IL}
    Let $X$, $Y$, and $Z$ be simply connected symplectic 
    $4$-manifolds with $b_{2}^{+}\equiv 3 \pmod 4$. If 
    $c_{1}^{2}(X)+c_{1}^{2}(Y)\leq 12$, then $X\# Y$ does not admit 
    Einstein metrics. Similarly, if 
    $c_{1}^{2}(X)+c_{1}^{2}(Y)+c_{1}^{2}(Z)\leq 24$, then $X\# Y\# Z$ 
    does not admit Einstein metrics.
    \end{cor}
\begin{proof}   
    By Taubes's result~\cite{T,Bourbaki}, the Chern classes $c_{1}(X)$ and 
    $c_{1}(Y)$ are basic classes with numerical Seiberg--Witten invariant 
    $=\pm 1$ on $X$ and $Y$ respectively. Applying Theorem~\ref{t:B} above, 
    we see that $c_{1}(X)+c_{1}(Y)$ is a monopole class on $X\# Y$, for which 
    the expected dimension of the monopole moduli space is $=1$. Therefore, 
    if $X\# Y$ admits an Einstein metric, Theorem~\ref{t:main} implies 
    $c_{1}(X)+c_{1}(Y)\geq 12$, with equality only if $X\# Y$ is symplectic, 
    which is not possible by~\cite{T,Bourbaki}.
    
    The proof of the second statement is entirely similar, but now the 
    expected dimension $d=2$.
\end{proof}
While this Corollary applies to spin manifolds, it is sometimes empty for 
non-obvious reasons.
\begin{rem}\label{r:homotopy}
    In the spin case the first part of Corollary~\ref{c:IL} is a consequence 
    of the Hitchin--Thorpe inequality~\eqref{eq:HT}, and is thus 
    homotopy-invariant. A spin symplectic $4$-manifold with 
    $b_{2}^{+}\equiv 3 \pmod 4$ must have $c_{1}^{2}\equiv 0\pmod{16}$ by 
    Rochlin's theorem. On the other hand, it is minimal, so that 
    $c_{1}^{2}\geq 0$ by the result of Taubes~\cite{Taubes,Bourbaki}. Therefore, 
    the assumption $c_{1}^{2}(X)+c_{1}^{2}(Y)\leq 12$ implies that $c_{1}^{2}$ 
    vanishes for both $X$ and $Y$, and then their connected sum violates the 
    Hitchin--Thorpe inequality~\eqref{eq:HT}.
\end{rem}
The second part of Corollary~\ref{c:IL} can be used for spin manifolds in 
exactly one way: one of the summands has $c_{1}^{2}=16$, and the other 
two have $c_{1}^{2}=0$.
\begin{ex}\label{kK3}
    Let $X$ be a symplectic spin manifold with $c_{1}^{2}(X)=16$ and 
    $\chi(X)=4$, where $\chi=\frac{1}{4}(e+\sigma)$ denotes the 
    holomorphic Euler characteristic. 
    Such manifolds exist by the results of Park and Szab\'o~\cite{PS}. By 
    Freedman's classification~\cite{freed}, such an $X$ is homeomorphic to 
    $K3\# 4(S^{2}\times S^{2})$. Take $Y=E(2n)$ a spin elliptic surface, and $Z$ 
    the symplectic spin manifold obtained from $E(2m)$ by performing a logarithmic 
    transformation of odd multiplicity $p$. Then Corollary~\ref{c:IL} applies to 
    show that the connected sum $X\# Y\# Z$ does not support an Einstein metric. 
    Note that $X\# Y\# Z$ is homeomorphic to 
    $(n+m+1)K3\# (n+m+2)(S^{2}\times S^{2})$. 
    As we increase $p$, the multiplicity of the logarithmic transformation, we 
    find that there are more and more basic classes on 
    $Z$ whose numerical Seiberg--Witten invariants are $\pm 1$, see 
    Fintushel--Stern~\cite{FSrat}, Theorem~8.7. By Theorem~\ref{t:B} these 
    basic classes give rise to special monopole classes on $X\# Y\# Z$, 
    so we have a sequence of smooth structures for which the number of 
    special monopole classes is unbounded. Thus Lemma~\ref{l:finite} 
    shows that we have infinitely many smooth structures.
    
    Varying $n$ and $m$, we see that $kK3\# (k+1)(S^{2}\times S^{2})$ has 
    infinitely many smooth structures not supporting Einstein metrics for 
    every $k\geq 3$. 
    \end{ex}
For $k\geq 4$ such examples were previously given by 
Ishida--LeBrun\footnote{The argument of Ishida--LeBrun~\cite{IL} for 
distinguishing infinitely many smooth structures is not correct as written, 
because the ``bandwidth'' is ill-defined taking the maximum over an infinite 
set, see Remark~\ref{r:IL} above. The argument can be corrected by considering 
special monopole classes only.}~\cite{IL} using building blocks due to 
Gompf. Whenever $k\geq 4$, we have no way of finding a smooth 
structure on $kK3\# (k+1)(S^{2}\times S^{2})$ supporting an Einstein metric 
because this manifold violates the Noether inequality (or the Miyaoka--Yau 
inequality, if one reverses the orientation). The case $k=3$ allows us to 
prove the following:
\begin{thm}\label{t:spin}
    There is a simply connected topological spin manifold $M$ which 
    satisfies the strict Hitchin--Thorpe inequality 
    $e (M)>\frac{3}{2}\vert\sigma (M)\vert$ and which admits 
    a smooth structure supporting an Einstein metric and also admits 
    infinitely many smooth structures not supporting Einstein metrics.
    \end{thm}
\begin{proof}
    We consider the manifold $M=3K3\# 4(S^{2}\times S^{2})$. We showed 
    in the above example that $M$ admits infinitely many smooth 
    structures without Einstein metrics. On the other hand, the 
    algebraic surface $S$ obtained as the double cover of the projective 
    plane branched in a smooth holomorphic curve of degree $10$ has 
    $c_{1}^{2}(S)=8$ and $\chi(S)=7$, and is spin because its canonical 
    bundle is the pullback of ${\mathcal O}(2)$. Thus $S$ is 
    homeomorphic to $M$ by~\cite{freed}. As it has ample canonical 
    bundle, it admits an Einstein metric by the results of 
    Aubin~\cite{A} and Yau~\cite{Y}. 
    \end{proof}
\begin{rem}
    The manifold $M$ has another infinite sequence of smooth structures, 
    which are distinct from the ones discussed above. 
    Fintushel--Stern~\cite{FScusp} have shown 
    that one can perform cusp surgery on a torus in $S$ to construct 
    infinitely many distinct smooth structures which are irreducible and 
    non-complex, and are therefore distinct from the smooth structures we 
    consider. Whether they admit Einstein metrics is not known.
\end{rem}

\begin{ex}\label{kkK3}
    Continuing Example~\ref{kK3}, we can prove that the following 
    topological spin manifolds admit infinitely many smooth structures not 
    supporting Einstein metrics: $kK3\# k(S^{2}\times S^{2})$ for all odd 
    $k\geq 5$, and $kK3\# (k+4)(S^{2}\times S^{2})$ for all even $k\geq 4$.

    We take the connected sum of one or two copies of the Park--Szab\'o 
    manifold $X$ used in Example~\ref{kK3}, and form the connected sum with 
    three respectively two simply connected spin elliptic surfaces. Using 
    the last part of Theorem~\ref{t:B} and Theorem~\ref{t:main} with $d=3$ 
    shows that such a connected sum has no Einstein metric if it has 
    $b_{2}^{+}\equiv 4\pmod8$. This covers all the cases mentioned above. 
    As before, we can obtain infinitely many distinct smooth structures 
    detected by monopole classes by performing logarithmic transformations 
    on the elliptic surfaces.

    These manifolds do not admit almost complex structures, so there is 
    certainly no chance to find K\"ahler--Einstein metrics for some other 
    smooth structure.
\end{ex}

Theorem~\ref{t:main} or Corollary~\ref{c:IL} can be applied to non-spin manifolds 
to produce some examples not admitting Einstein metrics and which have small 
homology.
\begin{ex}
    For every $q\geq 26$, the manifold $6\C P^{2}\# q\overline{\C 
    P^{2}}$ has smooth structures not admitting Einstein metrics. To 
    see this we can use the symplectic manifolds with $b_{2}^{+}=3$ 
    constructed by D.~Park~\cite{DP}. For example, he constructs a 
    minimal symplectic $4$-manifold $X$ homeomorphic to $3\C P^{2}\# 
    13\overline{\C P^{2}}$. This has $c_{1}^{2}(X)=6$. Applying the first 
    part of Corollary~\ref{c:IL} to this $X$ and blowups $Y$ of $X$ gives 
    the result.
    \end{ex}
Note that these manifolds have rather smaller homology than similar 
examples in~\cite{IL,BK}.

More interesting than the above are examples which can be shown to be 
homeomorphic to Einstein manifolds. The smallest example where the 
second part of Corollary~\ref{c:IL} applies is the following.
\begin{ex}\label{ex:quintic}
    The manifold $9\C P^{2}\# 44\overline{\C P^{2}}$ has a smooth 
    structure with an Einstein metric, and infinitely many without. 
    The smooth structure admitting an Einstein metric is that 
    underlying a smooth quintic in $\C P^{3}$. To obtain smooth 
    structures without Einstein metrics we can form the connected sum 
    of two copies of D.~Park's examples~\cite{DP}, with 
    $b_{2}^{-}=12$ and $13$, and one copy of a $K3$ surface. Then 
    performing logarithmic transformations on the $K3$ summand gives 
    infinitely many distinct symplectic structures not admitting 
    Einstein metrics as in Example~\ref{kK3} above.
    \end{ex}
This example is smaller than the ones in~\cite{K,lebrun3,IL}. It is 
interesting to note that one does not need any of the sophisticated ingredients 
we have used. It is easy to see that there are infinitely many distinct smooth 
structures on $9\C P^{2}\# 34\overline{\C P^{2}}$ which support symplectic 
forms and are distinguished by the numerical Seiberg--Witten invariants, 
see for example~\cite{P3}. The $10$-fold blowups of those smooth structures 
are still distinct, and do not admit Einstein metrics by~\cite{lebrun,K}. 
Thus, neither Theorem~\ref{t:lebrun}, nor any input from the stable cohomotopy 
refinement of the Seiberg--Witten invariants is needed for 
Example~\ref{ex:quintic}. 


\section{Towards irreducibility of Einstein manifolds}\label{s:irred}

Except for the connected sums $\C P^{2}\# k\overline{\C P^{2}}$ with 
$k\in\{1,3,\ldots,8\}$, all known Einstein four-manifolds are essentially 
irreducible. In fact, until recently the only other examples were either 
K\"ahler--Einstein, or of constant sectional 
curvature. Then, after the first version of this paper had been completed, 
Anderson~\cite{An} constructed some new examples of Einstein four-manifolds.
The manifolds in his examples also admit Riemannian metrics of non-positive 
sectional curvature. Thus, by the Cartan--Hadamard theorem, the following 
Lemma applies to them, as it does to the space forms of non-positive sectional 
curvature. The only orientable space form of positive curvature is $S^{4}$.
\begin{lem}\label{l:space}
    Let $M$ be any manifold with universal covering homeomorphic to 
    $\R^{4}$. Then $M$ is irreducible.
    \end{lem}
\begin{proof}
    Suppose $M=M_{1}\# M_{2}$. The fundamental group of $M$ 
    has only one end, and is therefore indecomposable as a free 
    product. Thus one of the $M_{i}$, say $M_{1}$, is simply 
    connected. Therefore $\pi_{1}(M_{2})=\pi_{1}(M)$, and so $M_{2}$ 
    carries the homology of $\pi_{1}(M)$. It follows that $M_{1}$ is a 
    homotopy sphere.
    \end{proof}

Concerning the K\"ahler--Einstein case, I proved in~\cite{Bourbaki} that 
if a minimal closed symplectic four-manifold $X$ with $b_{2}^{+}(X)>1$ 
decomposes as a connected sum, then one of the summands is an integral 
homology sphere whose fundamental group has no non-trivial finite quotients. 
In particular, if $\pi_{1}(X)$ is residually finite, then $X$ is 
irreducible. This applies to K\"ahler--Einstein surfaces of non-positive 
scalar curvature, because such manifolds are always minimal. Moreover, 
Gromov~\cite{gromov} proved that the fundamental group of a compact 
K\"ahler manifold never splits as a non-trivial free product. Thus, if 
one were to split off a non-trivial homology sphere from a K\"ahler 
surface, then the other summand would have to be simply connected.

As further, admittedly rather weak, evidence for the possible irreducibility 
of Einstein manifolds we have the following result.
\begin{thm}\label{t:nonE}
    Let $X$ be a smooth four-manifold with a monopole class 
    $c$. If $X$ admits an Einstein metric, then the maximal number $k$ of 
    copies of $\overline{\C P^{2}}$ that can be split off smoothly is 
    bounded by
    \begin{equation}\label{eq:main}
    k\leq\frac{1}{2}\left(2e(X)+3\sigma(X)-8d\right) \ .
    \end{equation}
    \end{thm}
\begin{proof}
    Suppose that $X\cong Y\#k\overline{\C P^{2}}$, and write 
    $c=c_{Y}+\sum_{i=1}^{k}a_{i}e_{i}$, with respect to the obvious 
    direct sum decomposition of $H^{2}(X,\Z)$. Here $e_{i}$ are the 
    generators for the cohomology of the $\overline{\C P^{2}}$ summands. 
    Note that the $a_{i}$ are odd integers because $c$ must be 
    characteristic. Now the reflections in the $e_{i}$ are realised by 
    self-diffeomorphisms of $X$, and the images of our monopole class under 
    these diffeomorphisms are again monopole classes. Thus, moving $c$ 
    by a diffeomorphism, we can arbitrarily change $e_{i}$ to its negative.
    
    Given an Einstein metric $g$ on $X$, we choose the signs in such a 
    way that $a_{i}e_{i}^{+}\cdot c_{Y}^{+}\geq 0$. Then we find
    \begin{align*}
    (c^{+})^{2} &= \left( c_{Y}^{+}+\sum_{i=1}^{k}a_{i}e_{i}^{+}\right)^{2}\geq 
    (c_{Y}^{+})^{2}\\
    &\geq c_{Y}^{2}=c^{2}+\sum_{i=1}^{k}a_{i}^{2} = 
    4d+2e(X)+3\sigma(X)+\sum_{i=1}^{k}a_{i}^{2} \ .
    \end{align*}
    On the other hand, applying~\eqref{eq:sharp} to $c$ and $g$ gives
    $$
    (c^{+})^{2}\leq\frac{3}{2}(2e(X)+3\sigma(X)) \ .
    $$
    Combining the two inequalities and noting that $a_{i}^{2}\geq 1$ 
    because all the $a_{i}$ are odd integers proves the result.
    \end{proof}
Parts of this argument are reminiscent of the proof of irreducibility of 
minimal symplectic manifolds with residually finite fundamental groups 
in~\cite{Bourbaki}, which however uses some much deeper 
ingredients~\cite{Taubes}. Note that one obtains a stronger inequality 
whenever not all the $a_{i}$ are $\pm 1$.


A precursor of the argument in the proof of Theorem~\ref{t:nonE} was first 
applied by LeBrun~\cite{lebrun} in the case where $X$ is the blowup of a 
K\"ahler or symplectic manifold $Y$, and $c$ is the basic class of the 
symplectic structure. I then observed in~\cite{K} that the argument works 
for blowups of arbitrary manifolds $Y$ with a basic class $c_{Y}$, because 
the blowup formula~\cite{FS,KMT} shows that $c_{Y}+\sum_{i=1}^{k}\pm e_{i}$ 
is a basic class on $X$. The formulation above is such that we circumvent 
the absence of a blowup formula for monopole classes.

\section*{Appendix: Ricci-flat four-manifolds}\label{s:ric}

\begin{center}
    by D.~Kotschick and J.~Wehrheim
\end{center}
\medskip

A manifold of dimension at least three that admits metrics of positive scalar 
curvature also admits metrics of vanishing scalar curvature. However, there 
are also manifolds which admit scalar-flat metrics although they do not have 
any of positive scalar curvature. Such manifolds are called strongly 
scalar-flat. The simplest examples are tori and other flat manifolds.

Bourguignon proved that scalar-flat metrics on strongly scalar-flat manifolds 
are in fact Ricci-flat, and therefore Einstein, compare~\cite{einstein}. In 
dimension three this implies that they are flat. In higher dimensions not 
much is known about strongly scalar-flat manifolds, although Futaki~\cite{Fut} 
and Dessai~\cite{Des} have proved some partial classification results in 
dimensions $>4$. 

In dimension $4$, the only known Ricci-flat manifolds are strongly scalar-flat. 
They are flat manifolds and finite quotients of $K3$ surfaces with Calabi--Yau 
metrics. The isometric quotients of $K3$ surfaces were classified by 
Hitchin~\cite{HT} in his discussion of the borderline case of the Hitchin--Thorpe 
inequality~\eqref{eq:HT}. He showed that the possible covering groups are 
$\Z_{2}$ and $\Z_{2}\times\Z_{2}$, both of which do actually occur. After 
earlier, unpublished, work of Calabi, Charlap--Sah and Levine, the closed 
orientable flat four-manifolds were classified by Hillman~\cite{Hil} and by 
Wagner~\cite{wagner}, who showed that there are $27$ distinct ones. By 
Bieberbach, all these manifolds are finite quotients of $T^{4}$. 

To investigate Ricci-flat four-manifolds, we first show that if there is a 
special monopole class on such a manifold, then it is one of the standard 
examples. The argument is a generalization of Witten's vanishing theorem for 
non-torsion monopole classes in the case of positive scalar curvature. 
\begin{prop}\label{t:Ric1}
    Let $X$ be a closed oriented $4$-manifold with a special monopole 
    class $c$. If $X$ admits a Ricci-flat metric $g$, then $(X,g)$ 
    is isometric to a finite quotient of $T^{4}$ or $K3$ with a 
    standard metric, and $c$ is a torsion class.
    \end{prop}
\begin{proof}
    As every solution of the monopole equations~\eqref{eq:D} 
    and~\eqref{eq:curv} on $(X,g)$ must be reducible, we conclude
    $c^{2}\leq 0$.
    
    On the other hand, $(X,g)$ is Einstein. Thus $X$ satisfies the 
    Hitchin--Thorpe inequality $2e(X)+3\sigma(X)\geq 0$, and the 
    assumption that $c$ is special gives $c^{2}\geq 2e(X)+3\sigma(X)\geq 0$.
    
    Combining the two inequalities, we conclude $c^{2}=0$. As the curvature 
    form $F_{\hat A}$ is anti-self-dual for every reducible solution 
    $(A,\Phi)$ of the monopole equations, we conclude that $F_{\hat A}$
    vanishes and $c$ is a torsion class.
    
    Moreover, we have $2e(X)+3\sigma(X) = 0$, so we are in the 
    borderline case of the Hitchin--Thorpe inequality~\eqref{eq:HT}, 
    and Hitchin's classification~\cite{HT} shows that every Einstein 
    metric on $X$ is either flat or an isometric quotient of a 
    Calabi--Yau metric on $K3$. 
    \end{proof}
This implies the following restrictions on any other Ricci-flat 
four-manifolds.
\begin{thm}\label{t:Ric2}
    If a closed orientable four-manifold $X$ admits a Ricci-flat metric 
    and is not a finite quotient of $T^{4}$ or $K3$, then:
    \begin{enumerate}
	\item it satisfies the strict Hitchin--Thorpe inequality 
	$2e(X)>3\vert\sigma(X)\vert$,
	\item it has finite fundamental group,
	\item $X$ and all its coverings have no special monopole classes for 
	either choice of orientation, and
	\item if the universal covering $\tilde X$ is spin, then $\sigma(X)=0$.
	\end{enumerate}
    \end{thm}
\begin{proof}
    We assume that $X$ admits an Einstein metric. As we have excluded the 
    manifolds occuring in the borderline case of the Hitchin--Thorpe 
    inequality~\eqref{eq:HT}, the inequality must be strict for $X$.
    In particular, the Euler characteristic of $X$ is positive. This, 
    together with the splitting theorem of Cheeger--Gromoll for 
    Ricci-flat manifolds, implies that $X$ has finite fundamental group.
  
    By Proposition~\ref{t:Ric1}, the assumptions imply that, for either 
    orientation, $X$ has no special monopole classes. The argument 
    also applies to any finite covering.
    
    Finally, suppose that $\tilde X$ is spin with non-zero signature. Then, 
    for any metric, there must be a not identically zero harmonic spinor of 
    pure  chirality. For a scalar-flat metric, the Weitzenb\"ock formula shows 
    that such a spinor $\Phi$ is parallel. Now choose the orientation 
    such that 
    $\Phi$ has positive chirality. Then $\sigma(\Phi,\Phi)$ is a non-zero 
    parallel self-dual two-form, i.~e.~a K\"ahler form. The K\"ahler 
    structure gives rise to a basic class $c_{1}(\tilde X)$. 
    
    Now $b_{1}(\tilde X)=0$, and $\vert\sigma(\tilde X)\vert\geq 16$ 
    by Rochlin's theorem, so that
    $$
    b_{2}^{+}(\tilde X)=\frac{1}{2}(e(\tilde X)+\sigma(\tilde 
    X))-1 > \frac{1}{4}\vert\sigma(\tilde X)\vert-1\geq 3
    $$
    using the strict Hitchin--Thorpe inequality for $\tilde X$. Thus 
    the basic class $c_{1}(\tilde X)$ is a special monopole class, 
    contradicting what we proved above.
     \end{proof}
The second conclusion of this theorem shows that a non-standard Ricci-flat 
four-manifold is not prevented from having positive scalar curvature by either 
enlargeability in the sense of Gromov--Lawson~\cite{GL}, or by the minimal 
hypersurface method of Schoen--Yau~\cite{SY}, because those only apply to 
manifolds with infinite fundamental groups. The last part shows that 
the Lichnerowicz obstruction vanishes as well. Note that the refined 
index-theoretic obstructions to positive scalar curvature due to 
Rosenberg reduce to the Lichnerowicz obstruction the case at hand, because 
the fundamental group is finite. As there are no special monopole classes 
for either orientation and all covering spaces, the only known obstruction 
to positive scalar curvature that could apply is the existence of a 
non-special monopole class. Here such a monopole class $c$ would actually 
have to satisfy $c^{2} < 0$, which is stronger than $c^{2} < 2e(X)+3\sigma(X)$ 
because $X$ satisfies the strict Hitchin--Thorpe inequality.
    
\bibliographystyle{amsplain}

\bigskip

\end{document}